\documentclass[10pt]{amsart}
\usepackage{amscd}
\usepackage{amssymb}

\newtheorem{thm}[equation]{Theorem}
\newtheorem{cor}[equation]{Corollary}
\newtheorem{prop}[equation]{Proposition}
\newtheorem{lem}[equation]{Lemma}
\theoremstyle{definition}
\newtheorem{dfn}[equation]{Definition}
\newtheorem{rem}[equation]{Remark}

\newtheorem{que}[equation]{Question}
\numberwithin{equation}{section}

\newcommand{\surj}{\twoheadrightarrow}

\newcommand{\ar}{\rightarrow}
\newcommand{\opn}{\operatorname}
\newcommand{\cat}[1]{\operatorname{\mathsf{#1}}}
\newcommand{\bdot}{{\textstyle \cdot}}

\newcommand{\rmitem}[1]{\item[\text{\textup{(#1)}}]}

\newcommand{\mfrak}[1]{\mathfrak{#1}}
\newcommand{\mcal}[1]{\mathcal{#1}}
\newcommand{\msf}[1]{\mathsf{#1}}
\newcommand{\mbf}[1]{\mathbf{#1}}
\newcommand{\mrm}[1]{\mathrm{#1}}
\newcommand{\mbb}[1]{\mathbb{#1}}

\newcommand{\tup}[1]{\textup{#1}}
\newcommand{\bsym}[1]{\boldsymbol{#1}}
\newcommand{\bwedge}{\bigwedge\nolimits}
\newcommand{\boplus}{\bigoplus\nolimits}
\newcommand{\what}{\widehat}

\title[Hochschild Cochain Complex]
{The Continuous Hochschild Cochain Complex of a Scheme}
\author{Amnon Yekutieli}
\address{Department of Mathematics,
Ben Gurion University, Be'er Sheva 84105, ISRAEL}
\email{amyekut@math.bgu.ac.il \newline \indent
http://www.math.bgu.ac.il/$\sim$amyekut}
\date{8.11.01}

\begin{document}

\begin{abstract}
Let $X$ be a separated finite type scheme over a noetherian base
ring  $\mbb{K}$. There is a complex $\what{\mcal{C}}^{\bdot}(X)$
of topological $\mcal{O}_{X}$-modules, called the complete
Hochschild chain complex of $X$. To any $\mcal{O}_{X}$-module
$\mcal{M}$ -- not necessarily quasi-coherent -- we assign the complex
$\mcal{H}om^{\mrm{cont}}_{\mcal{O}_{X}}
\bigl( \what{\mcal{C}}^{\bdot}(X), \mcal{M} \bigr)$
of continuous Hochschild cochains with values in $\mcal{M}$.
Our first main result is that when $X$ is smooth over $\mbb{K}$
there is a functorial isomorphism
\[ \mcal{H}om^{\mrm{cont}}_{\mcal{O}_{X}}
\bigl( \what{\mcal{C}}^{\bdot}(X), \mcal{M} \bigr)
\cong \mrm{R} \mcal{H}om_{\mcal{O}_{X^{2}}}
(\mcal{O}_{X}, \mcal{M}) \]
in the derived category
$\msf{D}(\cat{Mod} \mcal{O}_{X^{2}})$,
where $X^{2} := X \times_{\mbb{K}} X$.

The second main result is that if $X$ is smooth of relative
dimension $n$ and $n!$ is invertible in $\mbb{K}$, then
the standard maps
$\pi: \what{\mcal{C}}^{-q}(X) \to \Omega^{q}_{X / \mbb{K}}$
induce a quasi-isomorphism
\[ \mcal{H}om_{\mcal{O}_{X}}
\bigl( \boplus_{q} \Omega^{q}_{X / \mbb{K}}[q],
\mcal{M} \bigr)
\to \mcal{H}om^{\mrm{cont}}_{\mcal{O}_{X}}
\bigl( \what{\mcal{C}}^{\bdot}(X), \mcal{M} \bigr) . \]
When $\mcal{M} = \mcal{O}_{X}$ this is the quasi-isomorphism
underlying the Kontsevich Formality Theorem.

Combining the two results above we deduce a decomposition of the
global Hochschild cohomology
\[ \opn{Ext}^{i}_{\mcal{O}_{X^{2}}}
(\mcal{O}_{X}, \mcal{M}) \cong
\boplus_{q} \mrm{H}^{i - q}
\bigl( X, \bigl( \bwedge^{q}_{\mcal{O}_{X}} \mcal{T}_{X / \mbb{K}}
\bigr) \otimes_{\mcal{O}_{X}} \mcal{M} \bigr) , \]
where $\mcal{T}_{X / \mbb{K}}$ is the relative tangent sheaf.
\end{abstract}

\maketitle

\setcounter{section}{-1}
\section{Introduction and Statement of Results}

Let $\mbb{K}$ be a noetherian commutative ring and $X$ a
separated $\mbb{K}$-scheme of finite type. The diagonal morphism
$\Delta: X \to X^{2} = X \times_{\mbb{K}} X$ is then a closed
embedding. This allows us to identify the category
$\cat{Mod} \mcal{O}_{X}$ of $\mcal{O}_{X}$-modules with its
image inside $\cat{Mod} \mcal{O}_{X^{2}}$ under the functor
$\Delta_{*}$.

We shall use derived categories freely in this paper, following
the reference \cite{RD}.

\begin{dfn}[Hochschild Cohomology, First Definition] \label{dfn0.1}
\quad \quad \quad \quad
\begin{enumerate}
\item Given an $\mcal{O}_{X}$-module $\mcal{M}$ the Hochschild
cochain complex of $X$ with values in $\mcal{M}$ is
$\mrm{R} \mcal{H}om_{\mcal{O}_{X^{2}}}(\mcal{O}_{X}, \mcal{M})
\in \msf{D}(\cat{Mod} \mcal{O}_{X^{2}})$
\item The $q$th Hochschild cohomology of $X$ with values in
$\mcal{M}$ is
\[ \opn{Ext}^{q}_{\mcal{O}_{X^{2}}}(\mcal{O}_{X}, \mcal{M}) =
\mrm{H}^{q} \bigl( X^{2}, \mrm{R} \mcal{H}om_{\mcal{O}_{X^{2}}}
(\mcal{O}_{X}, \mcal{M}) \bigr) . \]
\end{enumerate}
\end{dfn}

This definition of Hochschild cohomology was considered by
Kontsevich \cite{Ko} and Swan \cite{Sw} among others. We observe
that if $\mbb{K}$ is a field, $A$ is a commutative $\mbb{K}$-algebra,
$A^{\mrm{e}} := A \otimes_{\mbb{K}} A$, $X := \opn{Spec} A$,
$M$ is an $A$-module and $\mcal{M}$ is the
quasi-coherent $\mcal{O}_{X}$-module associated to $M$, then
$\opn{Ext}^{q}_{\mcal{O}_{X^{2}}}(\mcal{O}_{X}, \mcal{M}) \cong
\opn{Ext}^{q}_{A^{\mrm{e}}}(A, M) \cong \mrm{HH}^{q}(A, M)$
is the usual Hochschild cohomology. This partly justifies the
definition. As we shall see, Definition \ref{dfn0.1} agrees with two
other plausible definitions of Hochschild cohomology of a scheme.

In Section \ref{sec1} we introduce the complex
$\what{\mcal{C}}^{\bdot}(X)$
of complete Hochschild chains of $X$. For any $q$ the sheaf
$\what{\mcal{C}}^{-q}(X) = \what{\mcal{C}}_{q}(X)$
is a topological $\mcal{O}_{X}$-module. (Note the unusual indexing,
due to our use of derived categories.)
If $q < 0$ then $\what{\mcal{C}}_{q}(X) = 0$, whereas for
any $q \geq 0$ and any affine open set
$U = \opn{Spec} A \subset X$
the group of sections $\Gamma(U, \what{\mcal{C}}_{q}(X))$
is an adic completion of the usual module of Hochschild chains
$\mcal{C}_{q}(A) = A^{\otimes (q + 2)} \otimes_{A^{\mrm{e}}} A$.
The coboundary operator
$\partial: \what{\mcal{C}}^{-q}(X) \ar \what{\mcal{C}}^{-q + 1}(X)$
is continuous.

\begin{dfn}[Hochschild Cohomology, Second Definition] \label{dfn0.2}
\quad \quad
\begin{enumerate}
\item Given an $\mcal{O}_{X}$-module $\mcal{M}$ the
continuous Hochschild cochain complex of $X$ with values in
$\mcal{M}$ is
$\mcal{H}om^{\mrm{cont}}_{\mcal{O}_{X}}
(\what{\mcal{C}}^{\bdot}(X), \mcal{M})$,
where $\mcal{M}$ has the discrete topology.
\item In the special case
$\mcal{M} = \mcal{O}_{X}$ we write
\[ \mcal{C}^{q}_{\mrm{cd}}(X) :=
\mcal{H}om^{\mrm{cont}}_{\mcal{O}_{X}}
(\what{\mcal{C}}_{q}(X), \mcal{O}_{X}) . \]
\item The $q$th Hochschild cohomology of $X$ with values in
$\mcal{M}$ is
\[ \mrm{H}^{q} \bigl( X, \mcal{H}om^{\mrm{cont}}_{\mcal{O}_{X}}
(\what{\mcal{C}}^{\bdot}(X), \mcal{M}) \bigr) . \]
\end{enumerate}
\end{dfn}

It turns out that on any open set $U$ as above we get
\[ \begin{aligned}
{} & \Gamma(U, \mcal{C}^{q}_{\mrm{cd}}(X)) \cong \\
{} & \quad \quad \{ f \in \opn{Hom}_{\mbb{K}}(A^{\otimes q}, A) \mid
f \text{ is a differential operator in each factor} \} .
\end{aligned} \]
Hence this is the same kind of Hochschild cochain complex considered
by Kontsevich in \cite{Ko}.

\begin{thm} \label{thm0.1}
Suppose $\mbb{K}$ is a noetherian ring and $X$ is a smooth
separated $\mbb{K}$-scheme. Given an $\mcal{O}_{X}$-module
$\mcal{M}$ there is an isomorphism
\[ \mcal{H}om^{\mrm{cont}}_{\mcal{O}_{X}}
\bigl( \what{\mcal{C}}^{\bdot}(X), \mcal{M} \bigr) \cong
\mrm{R} \mcal{H}om_{\mcal{O}_{X^{2}}}
(\mcal{O}_{X}, \mcal{M}) \]
in $\msf{D}(\cat{Mod} \mcal{O}_{X^{2}})$. This isomorphism is
functorial in $\mcal{M}$. In particular for
$\mcal{M} = \mcal{O}_{X}$ we get
\[ \mcal{C}^{\bdot}_{\mrm{cd}}(X) \cong
\mrm{R} \mcal{H}om_{\mcal{O}_{X^{2}}}
(\mcal{O}_{X}, \mcal{O}_{X}) . \]
\end{thm}

The theorem is proved in Section \ref{sec2}, where it is restated
as Corollary \ref{cor2.1}, and is deduced from the more general
Theorem \ref{thm2.1}.

Theorem \ref{thm0.1} says that on a smooth scheme the two
definitions of Hochschild cochain complexes coincide.
In Section \ref{sec3} we examine a third
definition of Hochschild cohomology, due to Swan \cite{Sw}. We
prove (Theorem \ref{thm3.4}) that when $X$ is flat over
$\mbb{K}$ this third definition also agrees with Definition
\ref{dfn0.1}.

In Section \ref{sec4} we look at the homomorphism
$\pi: \what{\mcal{C}}_{q}(X) \to \Omega^{q}_{X} =
\Omega^{q}_{X / \mbb{K}}$
given by the formula
\[ \pi((1 \otimes a_{1} \otimes \cdots \otimes a_{q} \otimes 1)
\otimes 1) = \mrm{d} a_{1} \wedge \cdots \wedge \mrm{d} a_{q} . \]
Let us denote by
$\mcal{T}_{X} = \mcal{T}_{X / \mbb{K}} :=
\mcal{H}om_{\mcal{O}_{X}}(\Omega^{1}_{X}, \mcal{O}_{X})$
the tangent sheaf, and
$\bwedge^{q} \mcal{T}_{X} :=
\bwedge_{\mcal{O}_{X}}^{q} \mcal{T}_{X}$.
Consider the complexes
$\boplus_{q} \Omega^{q}_{X / \mbb{K}}[q]$ and
$\boplus_{q} (\bwedge^{q} \mcal{T}_{X})[-q]$
with trivial co\-boundaries.

\begin{thm}[Decomposition] \label{thm0.2}
Let $\mbb{K}$ be a noetherian ring, let $X$ be a separated smooth
$\mbb{K}$-scheme of relative dimension $n$, and assume $n!$ is
invertible in $\mbb{K}$. Then for any
$\mcal{M} \in \cat{Mod} \mcal{O}_{X}$
the homomorphism of complexes
\[ \mcal{H}om_{\mcal{O}_{X}}
\bigl( \boplus_{q} \Omega^{q}_{X}[q], \mcal{M} \bigr)
\to \mcal{H}om^{\mrm{cont}}_{\mcal{O}_{X}}
\bigl( \what{\mcal{C}}^{\bdot}(X), \mcal{M} \bigr) \]
induced by $\pi$ is a quasi-isomorphism. In particular for
$\mcal{M} = \mcal{O}_{X}$ we get a quasi-isomorphism
\[ \pi_{\mrm{cd}}: \boplus_{q}
\bigl( \bwedge^{q} \mcal{T}_{X} \bigr)[-q]
\ar \mcal{C}^{\bdot}_{\mrm{cd}}(X) . \]
\end{thm}

Theorem \ref{thm0.2} is restated (in slightly more general form)
in Section \ref{sec4} as Theorem \ref{thm3.1} and proved there.

The quasi-isomorphism $\pi_{\mrm{cd}}$ underlies Kontsevich's
Formality Theorem. The fact that $\pi_{\mrm{cd}}$ is a
quasi-isomorphism in the case of a $C^{\infty}$ real manifold is
\cite{Ko} Theorem 4.6.11; cf.\ also \cite{Ts} Theorem 2.2.2.

Putting Theorems \ref{thm0.1} and \ref{thm0.2} together we obtain
a decomposition of the Hochs\-child cochain complex
\begin{equation} \label{eqn0.1}
\mrm{R} \mcal{H}om_{\mcal{O}_{X^{2}}}
(\mcal{O}_{X}, \mcal{M}) \cong
\boplus_{q} \bigl( \bwedge^{q} \mcal{T}_{X} \bigr)[-q]
\otimes_{\mcal{O}_{X}} \mcal{M}
\end{equation}
in $\msf{D}(\cat{Mod} \mcal{O}_{X^{2}})$.

Passing to global cohomology in (\ref{eqn0.1}) we obtain the
following corollary. It extends Corollary 2.6 of \cite{Sw} where
the assumptions are that $\mbb{K}$ is a field of characteristic
$0$ and $X$ is smooth and quasi-projective.

\begin{cor} \label{cor0.1}
Let $\mbb{K}$ be a noetherian ring, let $X$ be a separated smooth
$\mbb{K}$-scheme of relative dimension $n$, and assume $n!$ is
invertible in $\mbb{K}$. Then for any
$\mcal{M} \in \cat{Mod} \mcal{O}_{X}$ the Hochschild cohomology
decomposes:
\[ \opn{Ext}^{i}_{\mcal{O}_{X^{2}}}(\mcal{O}_{X}, \mcal{M})
\cong \boplus_{q} \mrm{H}^{i - q}
(X, \bigl( \bwedge^{q} \mcal{T}_{X} \bigr) \otimes_{\mcal{O}_{X}}
\mcal{M} \bigr) . \]
\end{cor}

Observe that for $\mcal{M} = \mcal{O}_{X}$,
$X$ affine and $A := \Gamma(X, \mcal{O}_{X})$
we recover the Hochschild-Kostant-Rosenberg Theorem
$\opn{Ext}^{i}_{A^{\mrm{e}}}(A, A) \cong
\bwedge^{i} \mcal{T}_{A}$.

\begin{rem}
This paper replaces ``Decomposition of the Hochschild Complex
of a Scheme in Arbitrary Characteristic,'' which has been withdrawn.
The proof of the main result of that paper, which relied on minimal
injective resolutions, turned out to have a serious gap in it. The
gap was discovered by M.\ Van den Bergh.
\end{rem}

\medskip \noindent \textbf{Acknowledgments.}\
The paper grew out of an inspiring series of  discussions with
Monique Lejeune-Jalabert in 1997 on the subject of characteristic
classes and Hochschild complexes. I wish to thank her, and the
Universit\'{e} Joseph Fourier, Grenoble, for their kind
hospitality. Also thanks to Vladimir Hinich, Colin Ingalls,
Joseph Lipman, Carlos Simpson and the referee for their helpful
suggestions. Finally I want to thank Michel Van den Bergh for
detecting an error in an earlier version of the
paper (see Remark above).

\section{Complete Hochschild Chains}
\label{sec1}

Let $\mbb{K}$ be a commutative ring and $A$ a commutative
$\mbb{K}$-algebra. As usual we write $A^{\mrm{e}} := A \otimes A$
where $\otimes := \otimes_{\mbb{K}}$. For any natural number $q$ let
$\mcal{B}_{q}(A) := A^{\otimes (q + 2)} =
A \otimes \cdots \otimes A$.
$\mcal{B}_{q}(A)$ is an $A^{\mrm{e}}$-module via the ring homomorphism
$a_{1} \otimes a_{2} \mapsto
a_{1} \otimes 1 \otimes \cdots \otimes 1 \otimes a_{2}$.
The (unnormalized) bar resolution is
\begin{equation} \label{eqn1.2}
\cdots \ar \mcal{B}_{2}(A) \xrightarrow{\partial}
\mcal{B}_{1}(A) \xrightarrow{\partial} \mcal{B}_{0}(A)
\ar A \ar 0 ,
\end{equation}
where $\partial$ is the $A^{\mrm{e}}$-linear homomorphism
\[ \partial(a_{0} \otimes \cdots \otimes a_{q + 1}) =
\sum_{i = 0}^{q} (-1)^{i} a_{0} \otimes \cdots \otimes
a_{i} a_{i + 1} \otimes \cdots \otimes a_{q + 1} . \]
The coboundary $\partial$ is $A^{\mrm{e}}$-linear, and the
complex (\ref{eqn1.2}) is split-exact with splitting homomorphism
$s(a_{0} \otimes \cdots \otimes a_{q + 1}) =
a_{0} \otimes \cdots \otimes a_{q + 1} \otimes 1$. The homomorphism
$s$ is $A$-linear  when $A$ acts via $a \mapsto a \otimes 1$.
Cf.\ \cite{Lo} Section 1.1.

For any $q$ let
\[ \mcal{C}_{q}(A) := \mcal{B}_{q}(A) \otimes_{A^{\mrm{e}}} A . \]
$\mcal{C}_{q}(A)$ is the module of degree $q$ Hochschild chains
of $A$.

Since we will be using derived categories, whose objects are
cochain complexes, we shall unfortunately have to abandon the
conventional notations for Hochschild chains. The first
departure will be to write the bar resolution as a cochain
complex, with $\mcal{B}^{-q}(A) := \mcal{B}_{q}(A)$. Likewise
we write $\mcal{C}^{-q}(A) := \mcal{C}_{q}(A)$.

  From now on $\mbb{K}$ is assumed to be a noetherian ring.
Let $A$ be a finitely generated $\mbb{K}$-algebra.
Denote by $I_{q}$ the kernel of the ring epimorphism
$\mcal{B}_{q}(A) \ar A$,
$a_{0} \otimes \cdots \otimes a_{q + 1} \mapsto
a_{0} \cdots a_{q + 1}$.
Let $\what{\mcal{B}}_{q}(A)$ be the $I_{q}$-adic completion of
$\mcal{B}_{q}(A)$. The homomorphisms $\partial$ and $s$ are
continuous for the $I_{q}$-adic topologies, and hence
$\what{\mcal{B}}^{\bdot}(A)$ is a complex and
$\what{\mcal{B}}^{\bdot}(A) \ar A$ is a continuously $A$-split
quasi-isomorphism. We call $\what{\mcal{B}}^{\bdot}(A)$ the
complete bar resolution.

Next define
\[ \what{\mcal{C}}^{-q}(A) = \what{\mcal{C}}_{q}(A) :=
\what{\mcal{B}}_{q}(A) \otimes_{\what{A^{\mrm{e}}}} A \cong
\what{\mcal{B}}_{q}(A) \otimes_{A^{\mrm{e}}} A  \]
and
\[ \mcal{C}^{q}_{\mrm{cd}}(A) :=
\opn{Hom}^{\mrm{cont}}_{\what{A^{\mrm{e}}}}(\what{\mcal{B}}_{q}(A)
, A) \cong
\opn{Hom}^{\mrm{cont}}_{A}(\what{\mcal{C}}_{q}(A), A)
, \]
where the superscript ``$\mrm{cont}$'' refers to continuous
homomorphisms with respect to the adic topology, and ``$\mrm{cd}$''
stands for ``continuous dual.''
We call $\what{\mcal{C}}_{q}(A)$ the module of complete
Hochschild chains, and $\mcal{C}^{q}_{\mrm{cd}}(A)$ the module of
continuous Hochschild cochains.

\begin{lem} \label{lem2.1}
Assume $A$ is flat over $\mbb{K}$. Then $\what{\mcal{B}}^{\bdot}(A)$
is a flat resolution of $A$ as $\what{\mcal{B}}_{0}(A)$-module.
\end{lem}

\begin{proof}
Let's write
$\mcal{B}_{q}(A) = \mcal{B}_{0}(A) \otimes A^{\otimes q}$.
Since $A$ is a flat $\mbb{K}$-algebra, it follows that
$\what{\mcal{B}}_{0}(A) \ar
\what{\mcal{B}}_{0}(A) \otimes A^{\otimes q}$
is flat . Now $\what{\mcal{B}}_{0}(A) \otimes A^{\otimes q}$
is noetherian, and $\what{\mcal{B}}_{q}(A)$ is an adic completion
of it, so
$\what{\mcal{B}}_{0}(A) \otimes A^{\otimes q} \ar
\what{\mcal{B}}_{q}(A)$
is flat.
\end{proof}

Suppose $Y$ is a noetherian scheme and $Y_{0} \subset Y$ is a
closed subset. The formal completion of $Y$ along $Y_{0}$ is a
noetherian formal scheme $\mfrak{Y}$ with underlying
topological space $Y_{0}$. The structure sheaf $\mcal{O}_{\mfrak{Y}}$
is a sheaf of topological rings with $\mcal{I}$-adic
topology, where $\mcal{I} \subset \mcal{O}_{Y}$ is any coherent
ideal sheaf defining the closed set $Y_{0}$.
The canonical morphism $\mfrak{Y} \to Y$ is flat, i.e.\
$\mcal{O}_{\mfrak{Y}}$ is a flat $\mcal{O}_{Y}$-algebra.
See \cite{EGA-I} Section 10.8 for details.

\begin{dfn} \label{dfn2.4}
Let $X$ be a finite type separated $\mbb{K}$-scheme.
For any $q \geq 2$ let $\mfrak{X}^{q}$ be the formal completion
of the scheme
$X^{q} := X \times_{\mbb{K}} \cdots \times_{\mbb{K}} X$
along the diagonal embedding of $X$.
\begin{enumerate}
\item For any $q \geq 0$ let
$\what{\mcal{B}}_{q}(X) := \mcal{O}_{\mfrak{X}^{q + 2}}$.
\item For any $q \geq 0$ the sheaf of degree $q$ complete
Hochschild chains of $X$ is
$\what{\mcal{C}}_{q}(X) :=
\what{\mcal{B}}_{ q}(X) \otimes_{\mcal{O}_{\mfrak{X}^{2}}}
\mcal{O}_{X}$.
\end{enumerate}
\end{dfn}

The benefit of the complete sheaves $\what{\mcal{B}}_{ q}(X)$ and
$\what{\mcal{C}}_{ q}(X)$ is they are coherent
(although over different ringed spaces). Indeed:

\begin{prop} \label{prop2.1}
On any affine open set $U = \opn{Spec} A \subset X$ one has
$\Gamma(U, \what{\mcal{B}}_{q}(X)) = \what{\mcal{B}}_{q}(A)$
and
$\Gamma(U, \what{\mcal{C}}_{q}(X)) = \what{\mcal{C}}_{q}(A)$.
\end{prop}

\begin{proof}
See \cite{EGA-I} Section 10.10.
\end{proof}

The homomorphisms
$\partial: \what{\mcal{B}}_{q}(A)\to \what{\mcal{B}}_{q - 1}(A)$
and
$s: \what{\mcal{B}}_{q}(A)\to \what{\mcal{B}}_{q + 1}(A)$
sheafify, hence
$\what{\mcal{B}}^{\bdot}(X)$ and $\what{\mcal{C}}^{\bdot}(X)$
are complexes with continuous coboundary operators, and
$\what{\mcal{B}}^{\bdot}(X) \ar \mcal{O}_{X}$ is a continuously
$\mcal{O}_{X}$-split quasi-isomorphism.

Given an $\mcal{O}_{X}$-module $\mcal{M}$ we consider
$\mcal{M} \otimes^{\mrm{L}}_{\mcal{O}_{X^{2}}} \mcal{O}_{X}
\cong \mrm{L} \Delta^{*} \mcal{M}$
as an object of $\msf{D}(\cat{Mod} \mcal{O}_{X})$.

\begin{prop} \label{prop2.4}
Assume $X$ is flat over $\mbb{K}$. Then there is an  isomorphism
\[ \what{\mcal{C}}^{\bdot}(X) \cong
\mcal{O}_{X} \otimes^{\mrm{L}}_{\mcal{O}_{X^{2}}} \mcal{O}_{X} \]
in $\msf{D}(\cat{Mod} \mcal{O}_{X})$.
\end{prop}

\begin{proof}
As in any completion of a noetherian scheme,
$\what{\mcal{B}}_{0}(X) = \mcal{O}_{\mfrak{X}^{2}}$ is a flat
$\mcal{O}_{X^{2}}$-algebra.
  From Lemma \ref{lem2.1} we see that $\what{\mcal{B}}_{q}(X)$
is a flat $\what{\mcal{B}}_{0}(X)$-module. Hence
$\what{\mcal{B}}^{\bdot}(X)$ is a flat resolution of
$\mcal{O}_{X}$ as $\mcal{O}_{X^{2}}$-module. But
$\what{\mcal{C}}^{\bdot}(X) \cong
\what{\mcal{B}}^{\bdot}(X) \otimes_{\mcal{O}_{X^{2}}}
\mcal{O}_{X}$.
\end{proof}

Given an $\mcal{O}_{X}$-module $\mcal{M}$ we have sheaves
$\mcal{H}om^{\mrm{cont}}_{\mcal{O}_{X}}
(\what{\mcal{C}}_{q}(X), \mcal{M})$,
where ``$\mcal{H}om^{\mrm{cont}}$'' refers to continuous
homomorphisms for the adic topology on $\what{\mcal{C}}_{q}(X)$
and the discrete topology on $\mcal{M}$. The continuous coboundary
$\partial$ makes
$\mcal{H}om^{\mrm{cont}}_{\mcal{O}_{X}}
(\what{\mcal{C}}^{\bdot}(X), \mcal{M})$
into a complex. In Definition \ref{dfn0.2} this was called the
continuous Hochschild cochain complex with values in
$\mcal{M}$.

\begin{prop} \label{prop2.2}
\begin{enumerate}
\item If $\mcal{M}$ is quasi-coherent then
$\mcal{H}om^{\mrm{cont}}_{\mcal{O}_{X}}
(\what{\mcal{C}}_{q}(X), \mcal{M})$
is also quasi-coherent.
\item For any affine open set
$U = \opn{Spec} A \subset X$, with
$M := \Gamma(U, \mcal{M})$, one has
\[ \Gamma \bigl( U, \mcal{H}om^{\mrm{cont}}_{\mcal{O}_{X}}
(\what{\mcal{C}}_{q}(X), \mcal{M}) \bigr) =
\opn{Hom}^{\mrm{cont}}_{A}
(\what{\mcal{C}}_{q}(A), M) . \]
\item With $U$ as above,
\[ \begin{aligned}
{} & \Gamma(U, \mcal{C}^{q}_{\mrm{cd}}(X)) =
\mcal{C}^{q}_{\mrm{cd}}(A) \cong \\
{} & \quad \quad \{ f \in \opn{Hom}_{\mbb{K}}(A^{\otimes q}, A) \mid
f \text{ is a differential operator in each factor} \} .
\end{aligned} \]
\end{enumerate}
\end{prop}

\begin{proof}
(1), (2) We have
\[ \mcal{H}om^{\mrm{cont}}_{\mcal{O}_{X}}
(\what{\mcal{C}}_{q}(X), \mcal{M}) \cong
\lim_{m \ar} \mcal{H}om_{\mcal{O}_{X^{2}}}
(\what{\mcal{B}}_{q}(X) / \what{\mcal{I}}_{q}^{m}, \mcal{M}) \]
where
$\what{\mcal{I}}_{q} = \opn{Ker}(\what{\mcal{B}}_{q}(X) \ar
\mcal{O}_{X})$.
But the sheaf $\what{\mcal{B}}_{q}(X) / \what{\mcal{I}}_{q}^{m}$
is a coherent $\mcal{O}_{X^{2}}$-module.

\medskip \noindent
(3) This is immediate from the results in \cite{EGA-IV}
Section 16.8.
\end{proof}

We see from part (3) of the proposition that this approach to
Hochschild cochains is the same as the one used by Kontsevich
\cite{Ko}.

\section{Comparison of Two Definitions}
\label{sec2}

In this section we prove that the two definitions of Hochschild
cohomology, Definitions \ref{dfn0.1} and \ref{dfn0.2}, coincide when
$X$ is smooth over $\mbb{K}$ (Corollary \ref{cor2.1}).
Throughout we assume $\mbb{K}$ is a noetherian ring and $X$
is a separated finite type scheme over $\mbb{K}$.

We start by recalling the notion of discrete
$\mcal{O}_{\mfrak{Y}}$-module on a noetherian formal scheme
$\mfrak{Y}$. An $\mcal{O}_{\mfrak{Y}}$-module $\mcal{M}$ is called
discrete if it is discrete for the adic topology of
$\mcal{O}_{\mfrak{Y}}$; in other words, if any local section of
$\mcal{M}$ is annihilated by some defining ideal of $\mfrak{Y}$.
The subcategory
$\cat{Mod}_{\mrm{disc}} \mcal{O}_{\mfrak{Y}} \subset
\cat{Mod} \mcal{O}_{\mfrak{Y}}$
of discrete modules is abelian and closed under direct limits.
Moreover $\cat{Mod}_{\mrm{disc}} \mcal{O}_{\mfrak{Y}}$ is locally
noetherian, so every injective object in
$\cat{Mod}_{\mrm{disc}} \mcal{O}_{\mfrak{Y}}$ is a direct sum of
indecomposable ones. The category
$\cat{Mod}_{\mrm{disc}} \mcal{O}_{\mfrak{Y}}$
has enough injectives, but we do not know if every injective in
$\cat{Mod}_{\mrm{disc}} \mcal{O}_{\mfrak{Y}}$
is also injective in the bigger category
$\cat{Mod}\mcal{O}_{\mfrak{Y}}$.
See \cite{RD} Section II.7 and \cite{Ye2} Sections 3-4 for details.

Given a point $y \in \mfrak{Y}$ let $\bsym{k}(y)$ be the residue
field and $\mcal{O}_{\mfrak{Y}, y}$ the local ring. Denote by
$\mcal{J}(y)$ an injective hull of $\bsym{k}(y)$ as
$\mcal{O}_{\mfrak{Y}, y}$-module. If $y'$ is a specialization of $y$
define $\mcal{J}(y, y')$ to be a constant sheaf on the
closed set $\overline{\{ y' \}}$ with stalk $\mcal{J}(y)$.

\begin{prop} \label{prop2.3}
Let $\mfrak{Y}$ be a noetherian formal scheme. The indecomposable
objects in $\cat{Mod}_{\mrm{disc}} \mcal{O}_{\mfrak{Y}}$ are the
sheaves $\mcal{J}(y, y')$.
\end{prop}

\begin{proof}
Exactly as in the proof of \cite{Ye2} Proposition 4.2.
\end{proof}

In particular this applies to $\mfrak{Y} = \mfrak{X}^{2}$, and we
shall denote by $\mcal{J}(x, x')$ the indecomposable
objects in $\cat{Mod}_{\mrm{disc}} \mcal{O}_{\mfrak{X}^{2}}$.
Therefore any injective $\mcal{J}$ in
$\cat{Mod}_{\mrm{disc}} \mcal{O}_{\mfrak{X}^{2}}$
has a decomposition
$\mcal{J} \cong \boplus_{x, x'} \mcal{J}(x, x')^{\mu(x, x')}$,
where $\mu(x, x')$ are cardinal numbers and
$\mcal{J}(x, x')^{\mu(x, x')}$ means a direct sum of $\mu(x, x')$
copies.

If $\mcal{M} \in \cat{Mod}_{\mrm{disc}} \mcal{O}_{\mfrak{X}^{2}}$
then
$\mcal{H}om^{\mrm{cont}}_{\mcal{O}_{\mfrak{X}^{2}}}
(\what{\mcal{B}}_{q}(X), \mcal{M})$
makes sense. The formula is
\[ \mcal{H}om^{\mrm{cont}}_{\mcal{O}_{\mfrak{X}^{2}}}
(\what{\mcal{B}}_{q}(X), \mcal{M}) =
\lim_{m \ar} \mcal{H}om_{\mcal{O}_{\mfrak{X}^{2}}}
\bigl( \what{\mcal{B}}_{q}(X) / \mcal{I}_{q}^{m},
\mcal{M} \bigr) \]
where
$\mcal{I}_{q} := \opn{Ker}(\what{\mcal{B}}_{q}(X) \ar
\mcal{O}_{X})$.
Hence given a complex
$\mcal{M}^{\bdot} \in \msf{D}(\cat{Mod}_{\mrm{disc}}
\mcal{O}_{\mfrak{X}^{2}})$
we obtain a total complex
$\mcal{H}om^{\mrm{cont}}_{\mcal{O}_{\mfrak{X}^{2}}}
(\what{\mcal{B}}^{\bdot}(X), \mcal{M}^{\bdot})$
with the usual indexing and signs.

Recall that $\mcal{O}_{\mfrak{X}^{2}}$ is an
$\mcal{O}_{X}$-algebra via the first projection
$X^{2} \ar X$, namely $a \mapsto a \otimes 1$.

\begin{lem} \label{lem2.4}
Let $\mcal{J}$ be an injective object
in $\cat{Mod}_{\mrm{disc}} \mcal{O}_{\mfrak{X}^{2}}$, and define
$\mcal{J}_{X} := \mcal{H}om_{\mcal{O}_{\mfrak{X}^{2}}}
(\mcal{O}_{X}, \mcal{J})$. Then there is a homomorphism of
$\mcal{O}_{X}$-modules $\tau: \mcal{J} \ar \mcal{J}_{X}$,
such that for any $q$ the induced homomorphism
\[ \tau_{q}: \mcal{H}om^{\mrm{cont}}_{\mcal{O}_{\mfrak{X}^{2}}}
(\what{\mcal{B}}_{q}(X), \mcal{J}) \ar
\mcal{H}om^{\mrm{cont}}_{\mcal{O}_{X}}
(\what{\mcal{B}}_{q}(X), \mcal{J}_{X}) \]
is an isomorphism.
\end{lem}

\begin{proof}
For any pair of points $x, x' \in X$ such that $x'$ is a
specialization of $x$ let
$\mcal{J}_{X}(x, x') \cong \mcal{H}om_{\mcal{O}_{\mfrak{X}^{2}}}
\bigl( \mcal{O}_{X}, \mcal{J}(x, x') \bigr)$
be the indecomposable injective $\mcal{O}_{X}$-module.

Let
$\mcal{I}_{q} := \opn{Ker}(\what{\mcal{B}}_{q}(X) \ar
\mcal{O}_{X})$,
a defining ideal of the formal scheme $\mfrak{X}^{q + 2}$.
For any $m \geq 1$ the sheaf of rings
$\what{\mcal{B}}_{q}(X) / \mcal{I}_{q}^{m}$
is coherent both as $\mcal{O}_{\mfrak{X}^{2}}$-module and as
$\mcal{O}_{X}$-module. We see that both
$\mcal{H}om_{\mcal{O}_{\mfrak{X}^{2}}}
\bigl(\what{\mcal{B}}_{q}(X)  / \mcal{I}_{q}^{m}, \mcal{J}(x, x')
\bigr)$ and \newline
$\mcal{H}om_{\mcal{O}_{X}}
( \what{\mcal{B}}_{q}(X) / \mcal{I}_{q}^{m},
\mcal{J}_{X}(x, x'))$
are constant sheaves on $\overline{\{ x' \}}$ with stalks being
injective hulls of $\bsym{k}(x)$ as
$(\what{\mcal{B}}_{q}(X) / \mcal{I}_{q}^{m})_{x}$-module.
Therefore
\[ \mcal{H}om_{\mcal{O}_{\mfrak{X}^{2}}}
\bigl( \what{\mcal{B}}_{q}(X) / \mcal{I}_{q}^{m}, \mcal{J}(x, x')
\bigr) \cong
\mcal{H}om_{\mcal{O}_{X}}
\bigl( \what{\mcal{B}}_{q}(X) / \mcal{I}_{q}^{m}, \mcal{J}_{X}(x, x')
\bigr) . \]
This isomorphism is not canonical, yet we can fit it into a
compatible direct system as $m$ varies. Thus there is a
(noncanonical) isomorphism
\begin{equation} \label{eqn2.1}
\begin{aligned}
\mcal{H}om^{\mrm{cont}}_{\mcal{O}_{\mfrak{X}^{2}}}
\bigl( \what{\mcal{B}}_{q}(X), \mcal{J}(x, x') \bigr) & \cong
\lim_{m \ar} \mcal{H}om_{\mcal{O}_{\mfrak{X}^{2}}}
\bigl( \what{\mcal{B}}_{q}(X) / \mcal{I}_{q}^{m},
\mcal{J}(x, x') \bigr) \\[1ex]
& \cong \lim_{m \ar} \mcal{H}om_{\mcal{O}_{X}}
\bigl( \what{\mcal{B}}_{q}(X) / \mcal{I}_{q}^{m},
\mcal{J}_{X}(x, x') \bigr) \\[1ex] & \cong
\mcal{H}om^{\mrm{cont}}_{\mcal{O}_{X}}
\bigl( \what{\mcal{B}}_{q}(X), \mcal{J}_{X}(x, x') \bigr) .
\end{aligned}
\end{equation}
Taking $q = 0$ above, and composing with homomorphism
``evaluation at $1$,'' we obtain
$\tau_{x, x'}: \mcal{J}(x, x') \ar \mcal{J}_{X}(x, x')$.

Now consider the given injective object $\mcal{J}$.
Choosing a decomposition
$\mcal{J} \cong \boplus_{x, x'} \mcal{J}(x, x')^{\mu(x, x')}$,
and summing up the homomorphisms $\tau_{x, x'}$,
we obtain a homomorphism $\tau: \mcal{J} \ar \mcal{J}_{X}$.
Because
\[ \mcal{H}om^{\mrm{cont}}_{\mcal{O}_{\mfrak{X}^{2}}}
(\what{\mcal{B}}_{q}(X), \mcal{J}) \cong
\boplus\nolimits_{x, x'}
\mcal{H}om^{\mrm{cont}}_{\mcal{O}_{\mfrak{X}^{2}}}
\bigl( \what{\mcal{B}}_{q}(X), \mcal{J}(x, x')
\bigr)^{\mu(x, x')} \]
and
\[ \mcal{H}om^{\mrm{cont}}_{\mcal{O}_{X}}
(\what{\mcal{B}}_{q}(X), \mcal{J}_{X}) \cong
\boplus\nolimits_{x, x'}
\mcal{H}om^{\mrm{cont}}_{\mcal{O}_{X}}
\bigl( \what{\mcal{B}}_{q}(X), \mcal{J}_{X}(x, x')
\bigr)^{\mu(x, x')} \]
it follows from (\ref{eqn2.1}) that $\tau_{q}$ is an isomorphism.
\end{proof}

Let $A$ be a $\mbb{K}$-algebra. For an element $a \in A$, an index
$q$ and any $1 \leq j \leq q$, let us define
\begin{equation} \label{eqn2.2}
\tilde{\mrm{d}}_{j} a :=
\underbrace{1 \otimes \cdots \otimes 1}_{j} \otimes
(a \otimes 1 - 1 \otimes a) \otimes 1 \otimes
\cdots \otimes 1 \in \mcal{B}_{q}(A) .
\end{equation}
Also let
\begin{equation} \label{eqn2.3}
\tilde{\mrm{d}}_{0} a :=
a \otimes 1 \otimes \cdots \otimes 1 -
1 \otimes \cdots \otimes 1 \otimes a \in \mcal{B}_{q}(A) .
\end{equation}
The ring $\mcal{B}_{q}(A)$ is an $A$-algebra by the map
$a \mapsto a \otimes 1 \otimes \cdots \otimes 1$.

Let $C$ be a noetherian commutative ring. A $C$-algebra $A$ is
\'etale if it is finitely generated and formally \'etale.

\begin{lem} \label{lem2.5}
Denote by $\mbb{K}[t_{1}, \ldots, t_{n}]$
the polynomial algebra in $n$ variables,
and let $\mbb{K}[t_{1}, \ldots, t_{n}] \to A$ be an \'{e}tale ring
homomorphism. Then for any $q \geq 0$ the ring
$\what{\mcal{B}}_{q}(A)$ is a formal power series algebra over
$A$ in the $n (q + 1)$ elements
$\tilde{\mrm{d}}_{j} t_{i}$.
\end{lem}

\begin{proof}
For any $q$ the homomorphism
$A \otimes \mbb{K}[t_{1}, \ldots, t_{n}]^{\otimes (q + 1)}
\to \mcal{B}_{q}(A)$
is \'{e}tale, which implies that
$A \to \mcal{B}_{q}(A)$ is formally smooth of relative
dimension $n (q + 1)$. In particular
$\Omega^{1}_{\mcal{B}_{q}(A) / A}$
is a free $\mcal{B}_{q}(A)$-module with basis
$\{ \mrm{d} t_{i, j} \}$, where $1 \leq j \leq q + 1$ and
\[ t_{i, j} :=
\underbrace{1 \otimes \cdots \otimes 1}_{j} \otimes t_{i}
\otimes 1 \otimes \cdots \otimes 1 \in \mcal{B}_{q}(A) . \]
Now for $1 \leq j \leq q$ we have
$\tilde{\mrm{d}}_{j} t_{i} = t_{i, j} - t_{i, j + 1}$, whereas
$\tilde{\mrm{d}}_{0} t_{i} - t_{i, q + 1} \in A$.
We see that the set
$\{ \mrm{d}(\tilde{\mrm{d}}_{j} t_{i}) \}_{j = 0}^{q}$
is also a basis of $\Omega^{1}_{\mcal{B}_{q}(A) / A}$.

Since the elements $\tilde{\mrm{d}}_{j} t_{i}$
are all in the defining ideal $I_{q}$ and since
$\mcal{B}_{q}(A) \ar \what{\mcal{B}}_{q}(A)$
is formally \'etale, we get a formally \'etale homomorphism
\[ \phi: A[[ \tilde{\mrm{d}}_{0} t_{1}, \ldots,
\tilde{\mrm{d}}_{q} t_{n} ]] \ar
\what{\mcal{B}}_{q}(A) . \]
Because $\phi$ lifts the identity $\phi_{0}: A \ar A$ it follows
that $\phi$ is bijective.
\end{proof}

Recall that $X$ is said to be smooth over $\mbb{K}$ if it is
formally smooth and finite type (see \cite{EGA-IV} Section 17). A
smooth scheme is also flat.

\begin{lem} \label{lem2.6}
Suppose $X$ is smooth over $\mbb{K}$. Then for any $q \geq 0$ the
functor
\[ \mcal{H}om^{\mrm{cont}}_{\mcal{O}_{\mfrak{X}^{2}}}
\bigl( \what{\mcal{B}}_{q}(X), - \bigr):
\cat{Mod}_{\mrm{disc}} \mcal{O}_{\mfrak{X}^{2}} \to
\cat{Mod}_{\mrm{disc}} \mcal{O}_{\mfrak{X}^{2}} \]
is exact.
\end{lem}

\begin{proof}
The statement can be verified locally on $X$, so let
$U = \opn{Spec} A \subset X$ be an affine open set that is \'{e}tale
over affine space $\mbf{A}^{n}_{\mbb{K}}$; cf.\ \cite{EGA-IV}
Corollary 17.11.3. In other
words there is an \'{e}tale ring homomorphism
$\mbb{K}[t_{1}, \ldots, t_{n}] \ar A$.
According to Lemma \ref{lem2.5}, $\what{\mcal{B}}_{q}(A)$ is a
formal power series algebra over
$\what{A^{\mrm{e}}} = \what{\mcal{B}}_{0}(A)$ in the elements
$\tilde{\mrm{d}}_{j} t_{i}$, where $1 \leq j \leq q$.

Denote by $I_{q, \mrm{e}}$ the kernel of the ring homomorphism
$\mcal{B}_{q}(A) \ar A^{\mrm{e}}$,
$a_{0} \otimes a_{1} \otimes \cdots \otimes a_{q + 1}
\mapsto a_{0} \otimes a_{1} \cdots a_{q + 1}$.
Let $\what{I}_{q, \mrm{e}}$ be its completion.
For any $m \geq 0$ the $\what{A^{\mrm{e}}}$-module
$\what{\mcal{B}}_{q}(A) / \what{I}_{q, \mrm{e}}^m$ is free of
finite rank -- with basis consisting of monomials in the
$\tilde{\mrm{d}}_{j} t_{i}$ -- and it
has the $\what{I}_{0}$-adic topology.

Passing to sheaves we see that for any $m$ the functor
$\mcal{H}om_{\mcal{O}_{\mfrak{X}^{2}}}
\bigl( \what{\mcal{B}}_{q}(X) /
\what{\mcal{I}}_{q, \mrm{e}}^m, - \bigr)$
is exact. But for any discrete module $\mcal{M}$,
\[ \mcal{H}om^{\mrm{cont}}_{\mcal{O}_{\mfrak{X}^{2}}}
\bigl( \what{\mcal{B}}_{q}(X) , \mcal{M} \bigr) \cong
\lim_{m \ar} \mcal{H}om_{\mcal{O}_{\mfrak{X}^{2}}}
\bigl( \what{\mcal{B}}_{q}(X) /
\what{\mcal{I}}_{q, \mrm{e}}^m, \mcal{M} \bigr) . \]
\end{proof}

\begin{thm} \label{thm2.1}
Suppose $\mbb{K}$ is a noetherian ring and $X$ is a smooth
separated $\mbb{K}$-scheme. Given a complex
$\mcal{M}^{\bdot} \in
\msf{D}^{+}(\cat{Mod}_{\mrm{disc}} \mcal{O}_{\mfrak{X}^{2}})$
there is an isomorphism
\[ \mcal{H}om^{\mrm{cont}}_{\mcal{O}_{\mfrak{X}^{2}}}
(\what{\mcal{B}}^{\bdot}(X), \mcal{M}^{\bdot})
\cong
\mrm{R} \mcal{H}om_{\mcal{O}_{X^{2}}}
(\mcal{O}_{X}, \mcal{M}^{\bdot}) \]
in $\msf{D}(\cat{Mod} \mcal{O}_{X^{2}})$. This isomorphism is
functorial in $\mcal{M}^{\bdot}$.
\end{thm}

\begin{proof}
Let $\mcal{M}^{\bdot} \ar \mcal{J}^{\bdot}$
be an injective resolution of $\mcal{M}^{\bdot}$ in
$\cat{Mod}_{\mrm{disc}} \mcal{O}_{\mfrak{X}^{2}}$.
By this we mean that $\mcal{M}^{\bdot} \ar \mcal{J}^{\bdot}$ is a
quasi-isomorphism and $\mcal{J}^{\bdot}$ is a bounded below
complex of injectives objects in
$\cat{Mod}_{\mrm{disc}} \mcal{O}_{\mfrak{X}^{2}}$.
Then each $\mcal{J}^{q}$ is an injective $\mcal{O}_{X^{2}}$-module
supported on $X$, and
\[ \mcal{H}om^{\mrm{cont}}_{\mcal{O}_{\mfrak{X}^{2}}}
(\mcal{O}_{X}, \mcal{J}^{\bdot}) =
\mcal{H}om_{\mcal{O}_{X^{2}}}(\mcal{O}_{X}, \mcal{J}^{\bdot})
= \mrm{R} \mcal{H}om_{\mcal{O}_{X^{2}}}
(\mcal{O}_{X}, \mcal{M}^{\bdot}) . \]

Since the homomorphism
$\what{\mcal{B}}^{\bdot}(X) \ar \mcal{O}_{X}$
is split by the continuous $\mcal{O}_{X}$-linear homomorphism
$s$, Lemma \ref{lem2.4} says that for any $q \geq 0$ the
homomorphism
\[ \mcal{H}om^{\mrm{cont}}_{\mcal{O}_{\mfrak{X}^{2}}}
(\mcal{O}_{X}, \mcal{J}^{q})
\ar \mcal{H}om^{\mrm{cont}}_{\mcal{O}_{\mfrak{X}^{2}}}
(\what{\mcal{B}}^{\bdot}(X), \mcal{J}^{q}) \]
is a quasi-isomorphism. Because $\what{\mcal{B}}^{\bdot}(X)$
is bounded above and $\mcal{J}^{\bdot}$ is bounded below, the
usual spectral sequence shows that
\[ \mcal{H}om^{\mrm{cont}}_{\mcal{O}_{\mfrak{X}^{2}}}
(\mcal{O}_{X}, \mcal{J}^{\bdot})
\ar \mcal{H}om^{\mrm{cont}}_{\mcal{O}_{\mfrak{X}^{2}}}
(\what{\mcal{B}}^{\bdot}(X), \mcal{J}^{\bdot}) \]
is a quasi-isomorphism.

Next by Lemma \ref{lem2.6} for any $q \leq 0$ the homomorphism
\[ \mcal{H}om^{\mrm{cont}}_{\mcal{O}_{\mfrak{X}^{2}}}
(\what{\mcal{B}}^{q}(X), \mcal{M}^{\bdot}) \ar
\mcal{H}om^{\mrm{cont}}_{\mcal{O}_{\mfrak{X}^{2}}}
(\what{\mcal{B}}^{q}(X), \mcal{J}^{\bdot}) \]
is a quasi-isomorphism. Therefore
\[ \mcal{H}om^{\mrm{cont}}_{\mcal{O}_{\mfrak{X}^{2}}}
(\what{\mcal{B}}^{\bdot}(X), \mcal{M}^{\bdot}) \ar
\mcal{H}om^{\mrm{cont}}_{\mcal{O}_{\mfrak{X}^{2}}}
(\what{\mcal{B}}^{\bdot}(X), \mcal{J}^{\bdot}) \]
is a quasi-isomorphism.
\end{proof}

Now we may compare the two definitions of Hochschild cochain
complexes.

\begin{cor} \label{cor2.1}
Suppose $\mbb{K}$ is a noetherian ring and $X$ is a smooth
separated $\mbb{K}$-scheme. Given an $\mcal{O}_{X}$-module
$\mcal{M}$ there is an isomorphism
\[ \mcal{H}om^{\mrm{cont}}_{\mcal{O}_{X}}
(\what{\mcal{C}}^{\bdot}(X), \mcal{M}) \cong
\mrm{R} \mcal{H}om_{\mcal{O}_{X^{2}}}
(\mcal{O}_{X}, \mcal{M}) \]
in $\msf{D}(\cat{Mod} \mcal{O}_{X^{2}})$. This isomorphism is
functorial in $\mcal{M}$. In particular for
$\mcal{M} = \mcal{O}_{X}$ we get
\[ \mcal{C}^{\bdot}_{\mrm{cd}}(X) \cong
\mrm{R} \mcal{H}om_{\mcal{O}_{X^{2}}}
(\mcal{O}_{X}, \mcal{O}_{X}) . \]
\end{cor}

\begin{proof}
This is immediate from the theorem, since
$\cat{Mod} \mcal{O}_{X} \subset
\cat{Mod}_{\mrm{disc}} \mcal{O}_{\mfrak{X}^{2}}$,
and
$\what{\mcal{C}}^{\bdot}(X) = \mcal{O}_{X} \otimes_{\mcal{O}_{X^{2}}}
\what{\mcal{B}}^{\bdot}(X)$.
\end{proof}

\begin{rem}
Assume $\mbb{K}$ is a field, and let $\mcal{K}^{\bdot}_{X}$ be
the residue complex of $X$, see \cite{Ye3}.
If $X$ is smooth of dimension $n$ over $\mbb{K}$ then
$0 \to \Omega^{n}_{X} \to \mcal{K}_{X}^{-n} \to \cdots \to
\mcal{K}_{X}^{0} \to 0$
is a minimal injective resolution. Hence
$\mcal{H}om^{\mrm{cont}}_{\mcal{O}_{X}}
(\what{\mcal{C}}^{\bdot}(X), \mcal{K}^{\bdot}_{X})$
is a bounded below complex of flasque sheaves isomorphic to
$\mrm{R} \mcal{H}om_{\mcal{O}_{X^{2}}}
(\mcal{O}_{X}, \Omega^{n}_{X})[n]$.
Moreover if $f: X \to Y$ is a proper morphism the trace
$\opn{Tr}_{f}: f_{*} \mcal{K}_{X} \to \mcal{K}_{Y}$
induces a homomorphism of complexes
\[ f_{*} \mcal{H}om^{\mrm{cont}}_{\mcal{O}_{X}}
(\what{\mcal{C}}^{\bdot}(X), \mcal{K}^{\bdot}_{X}) \to
\mcal{H}om^{\mrm{cont}}_{\mcal{O}_{Y}}
(\what{\mcal{C}}^{\bdot}(Y), \mcal{K}^{\bdot}_{Y}) . \]
This angle ought to be explored.
\end{rem}

\section{A Third Definition}
\label{sec3}

In the paper \cite{Sw} Swan makes the following definition. Let
$\mbb{K}$ be a commutative ring and $X$ a $\mbb{K}$-scheme. Let
$\mcal{C}_{q}(X)$ be the sheaf on $X$ associated to the presheaf
$U \mapsto \mcal{C}_{q} \bigl( \Gamma(U, \mcal{O}_{X}) \bigr)$.
Then $\mcal{C}^{\bdot}(X)$ is a complex of $\mcal{O}_{X}$-modules.
Given an $\mcal{O}_{X}$-module $\mcal{M}$ choose an injective
resolution
$\mcal{M} \ar \mcal{J}^{0} \ar \mcal{J}^{1} \ar \cdots$.
The $q$th Hochschild cohomology of $X$ with values in $\mcal{M}$ is
defined to be
\[ \mrm{H} \mrm{H}^{q}(X, \mcal{M}) :=
\mrm{H}^{q} \Gamma(X, \mcal{H}om_{\mcal{O}_{X}}
(\mcal{C}^{\bdot}(X), \mcal{J}^{\bdot})) . \]

This section is devoted to proving the following theorem.

\begin{thm} \label{thm3.4}
Let $\mbb{K}$ be a noetherian ring and $X$ a flat finite type
separated $\mbb{K}$-scheme. Let
$\mcal{M}^{\bdot} \in \msf{D}^{+}(\cat{Mod} \mcal{O}_{X})$
be a complex. Assume either of the following:
\begin{enumerate}
\rmitem{i} $X$ is embeddable as a closed subscheme of some smooth
$\mbb{K}$-scheme, and $\mbb{K}$ is a regular ring.
\rmitem{ii} Each $\mrm{H}^{q} \mcal{M}^{\bdot}$ is quasi-coherent.
\end{enumerate}
Then there is an isomorphism
\[ \mrm{R} \mcal{H}om_{\mcal{O}_{X}}
(\mcal{C}^{\bdot}(X), \mcal{M}^{\bdot}) \cong
\mrm{R} \mcal{H}om_{\mcal{O}_{X^{2}}}
(\mcal{O}_{X}, \mcal{M}^{\bdot}) \]
in $\msf{D}^{+}(\cat{Mod} \mcal{O}_{X^{2}})$. This isomorphism is
functorial in $\mcal{M}^{\bdot}$.
\end{thm}

\begin{cor} \label{cor3.4}
Under the assumptions of the theorem, with
$\mcal{M}^{\bdot} = \mcal{M}$ a single $\mcal{O}_{X}$-module,
there is an isomorphism
\[ \mrm{H} \mrm{H}^{q}(X, \mcal{M}) \cong
\opn{Ext}^{q}_{\mcal{O}_{X^{2}}}(\mcal{O}_{X}, \mcal{M}) . \]
\end{cor}

Corollary \ref{cor3.4} was proved by Swan in the case of a
field $\mbb{K}$ and a quasi-projective scheme $X$
(\cite{Sw} Theorem 2.1).

The proofs of Theorem \ref{thm3.4} and Corollary
\ref{cor3.4} are at the end of the section, after some
preparation.

The sheaves $\mcal{C}_{q}(X)$ are ill behaved; they are not
quasi-coherent except in trivial cases. The sheaves
$\mcal{B}_{q}(X)$, associated to the presheaves
$U \mapsto \mcal{B}_{q}(\Gamma(U, \mcal{O}_{X}))$,
are even more troublesome: we do not know if $\mcal{B}_{q}(X)$ is
an $\mcal{O}_{X^{2}}$-module. We get around these problems by
using the completions $\what{\mcal{C}}_{q}(X)$.

\begin{prop} \label{prop3.4}
Let $\mbb{K}$ be a noetherian ring and $X$ a flat finite type
separated $\mbb{K}$-scheme. Then there is an isomorphism
\[ \mcal{C}^{\bdot}(X) \cong \mcal{O}_{X}
\otimes^{\mrm{L}}_{\mcal{O}_{X^{2}}} \mcal{O}_{X} \]
in $\msf{D}(\cat{Mod} \mcal{O}_{X})$.
\end{prop}

\begin{proof}
For any affine open set $U = \opn{Spec} A \subset X$ there is a
quasi-isomorphism
$\mcal{C}^{\bdot}(A) \ar \what{\mcal{C}}^{\bdot}(A)$;
see Lemma \ref{lem2.1}.
Therefore when we pass to sheaves we obtain a quasi-isomorphism
$\mcal{C}^{\bdot}(X) \ar \what{\mcal{C}}^{\bdot}(X)$.
Now use Proposition \ref{prop2.4}.
\end{proof}

\begin{dfn} \label{dfn3.1}
Let $Y$ be a noetherian scheme.
An $\mcal{O}_{Y}$-module $\mcal{L}$ is called {\em finite pseudo
locally free} if
$\mcal{L} \cong \boplus_{i = 1}^{n} g_{i !} \mcal{L}_{i}$,
where for each $i$, $g_{i} : U_{i} \to Y$ is the inclusion of an
affine open set, $g_{i !}$ is extension by zero, and
$\mcal{L}_{i}$ is a locally free $\mcal{O}_{U_{i}}$-module of
finite rank.
\end{dfn}

According to \cite{RD} Theorem II.7.8, for any noetherian scheme
$Y$ the category $\cat{Mod} \mcal{O}_{Y}$ is locally noetherian.

\begin{lem} \label{lem3.4}
Suppose $Y$ is a noetherian scheme.
\begin{enumerate}
\item A finite pseudo locally free $\mcal{O}_{Y}$-module
$\mcal{L}$ is a noetherian object in
$\cat{Mod} \mcal{O}_{Y}$.
\item Given a noetherian object $\mcal{M} \in \cat{Mod} \mcal{O}_{Y}$
there is a surjection $\mcal{L} \surj \mcal{M}$ with $\mcal{L}$ a
finite pseudo locally free $\mcal{O}_{Y}$-module.
\item Let $\mcal{L}$ be a finite pseudo locally free
$\mcal{O}_{Y}$-module. Then $\mcal{L}$ is a flat
$\mcal{O}_{Y}$-module.
\item If $Y$ is separated and $\mcal{L}$ is a finite pseudo locally
free $\mcal{O}_{Y}$-module then the functor
\[ \mcal{H}om_{\mcal{O}_{Y}}(\mcal{L}, -):
\cat{QCoh} \mcal{O}_{Y} \to \cat{QCoh} \mcal{O}_{Y} \]
is exact.
\end{enumerate}
\end{lem}

\begin{proof}
(1) By the proof of \cite{RD} Theorem II.7.8, for any inclusion
$g: U \to X$ of an affine open set, the sheaf $g_{!} \mcal{O}_{U}$ is
noetherian in $\cat{Mod} \mcal{O}_{Y}$. This implies that for any
coherent $\mcal{O}_{U}$-module $\mcal{M}$, $g_{!} \mcal{M}$ is
noetherian in $\cat{Mod} \mcal{O}_{Y}$.

\medskip \noindent (2)
For every affine open subset $g: U \to Y$ and every section of
$\Gamma(U, \mcal{M})$ we get a homomorphism
$g_{!} \mcal{O}_{U} \to \mcal{M}$. By the ascending chain
condition finitely many of these cover $\mcal{M}$.

\medskip \noindent (3)
In order to verify flatness we may restrict to a
sufficiently small open subset $V \subset Y$. Thus
we can assume each $\mcal{L}_{i}$ in Definition \ref{dfn3.1}
is free; and hence we reduce to
the case $\mcal{L} = g_{!} \mcal{O}_{U}$ for an affine open subset
$g: U \to Y$.

For any $\mcal{O}_{Y}$-module $\mcal{M}$ we have
\[ g_{!} \mcal{O}_{U} \otimes_{\mcal{O}_{Y}} \mcal{M} \cong
g_{!} g^{*} \mcal{M} . \]
Since both functors $g^{*}$ and $g_{!}$ are exact it follows that
$g_{!} \mcal{O}_{U}$ is flat.

\medskip \noindent (4)
After the same reduction as in 3 we have
\[ \mcal{H}om_{\mcal{O}_{Y}}(g_{!} \mcal{O}_{U}, \mcal{M})
\cong g_{*} g^{*} \mcal{M} . \]
Since $g$ is now an affine morphism the functor
\[ g_{*}: \cat{QCoh} \mcal{O}_{U} \to \cat{QCoh} \mcal{O}_{Y} \]
is exact.
\end{proof}

\begin{proof}[Proof of Theorem \tup{\ref{thm3.4}}]
If condition (i) is satisfied then $X^{2}$ is embeddable in a
regular scheme. Hence we can find a resolution
$\cdots \to \mcal{L}^{-1} \to \mcal{L}^{0} \to \mcal{O}_{X}$
where all the $\mcal{O}_{X^{2}}$-modules $\mcal{L}^{q}$ are
locally free of finite rank. Otherwise by Lemma \ref{lem3.4}
we can at least find such a resolution where the $\mcal{L}^{q}$ are
finite pseudo locally free $\mcal{O}_{X^{2}}$-modules.

Since $\mcal{L}^{\bdot}$ is a flat resolution of $\mcal{O}_{X}$,
by Proposition \ref{prop3.4} we have
\[ \mcal{C}^{\bdot}(X) \cong
\mcal{O}_{X} \otimes_{\mcal{O}_{X^{2}}} \mcal{L}^{\bdot} \]
in $\msf{D}^{-}(\cat{Mod} \mcal{O}_{X})$.

Choose a quasi-isomorphism
$\mcal{M}^{\bdot} \to \mcal{K}^{\bdot}$
where $\mcal{K}^{\bdot}$ is a bounded below complex of
injective $\mcal{O}_{X}$-modules. Then choose a quasi-isomorphism
$\mcal{K}^{\bdot} \to \mcal{J}^{\bdot}$
where $\mcal{J}^{\bdot}$ is a bounded below complex of
injective $\mcal{O}_{X^{2}}$-modules. If condition (ii) holds then
take $\mcal{K}^{\bdot}$ and $\mcal{J}^{\bdot}$ to be complexes of
quasi-coherent injective modules over $\mcal{O}_{X}$ and
$\mcal{O}_{X^{2}}$ respectively (cf.\ \cite{RD} Theorem II.7.18).

We have
\[ \mrm{R} \mcal{H}om_{\mcal{O}_{X^{2}}}
(\mcal{O}_{X}, \mcal{M}^{\bdot}) =
\mcal{H}om_{\mcal{O}_{X^{2}}} (\mcal{O}_{X}, \mcal{J}^{\bdot}) , \]
and there is a quasi-isomorphism
\[ \mcal{H}om_{\mcal{O}_{X^{2}}} (\mcal{O}_{X}, \mcal{J}^{\bdot})
\to
\mcal{H}om_{\mcal{O}_{X^{2}}} (\mcal{L}^{\bdot}, \mcal{J}^{\bdot})
. \]
Since either all the $\mcal{L}^{q}$ are locally free
$\mcal{O}_{Y}$-modules of finite
rank (in case condition (i) holds), or all the $\mcal{L}^{q}$ are
finite pseudo locally free and all the $\mcal{K}^{q}$ and
$\mcal{J}^{q}$ are quasi-coherent (in case condition (ii) holds),
it follows that we have a quasi-isomorphism
\[ \mcal{H}om_{\mcal{O}_{X^{2}}} (\mcal{L}^{\bdot}, \mcal{K}^{\bdot})
\to
\mcal{H}om_{\mcal{O}_{X^{2}}} (\mcal{L}^{\bdot}, \mcal{J}^{\bdot})
. \]
But
\[ \mcal{H}om_{\mcal{O}_{X^{2}}} (\mcal{L}^{\bdot}, \mcal{K}^{\bdot})
\cong
\mcal{H}om_{\mcal{O}_{X}} (\mcal{O}_{X} \otimes_{\mcal{O}_{X^{2}}}
\mcal{L}^{\bdot}, \mcal{K}^{\bdot}) . \]
Finally
\[ \begin{aligned}
\mcal{H}om_{\mcal{O}_{X}} (\mcal{O}_{X} \otimes_{\mcal{O}_{X^{2}}}
\mcal{L}^{\bdot}, \mcal{K}^{\bdot}) & =
\mrm{R} \mcal{H}om_{\mcal{O}_{X}}
(\mcal{O}_{X} \otimes_{\mcal{O}_{X^{2}}}
\mcal{L}^{\bdot}, \mcal{K}^{\bdot}) \\
& \cong \mrm{R} \mcal{H}om_{\mcal{O}_{X}}
(\mcal{C}^{\bdot}(X), \mcal{K}^{\bdot}) \\
& = \mrm{R} \mcal{H}om_{\mcal{O}_{X}} (\mcal{C}^{\bdot}(X),
  \mcal{M}^{\bdot})
\end{aligned} \]
in $\msf{D}(\cat{Mod} \mcal{O}_{X})$. To this isomorphism we apply
the functor $\Delta_{*}$.
\end{proof}

\begin{proof}[Proof of Corollary \tup{\ref{cor3.4}}]
Choose an injective resolution $\mcal{M} \ar \mcal{J}^{\bdot}$.
Then
\[ \mcal{H}om_{\mcal{O}_{X}} (\mcal{C}^{\bdot}(X),
\mcal{J}^{\bdot}) =
\mrm{R} \mcal{H}om_{\mcal{O}_{X}} (\mcal{C}^{\bdot}(X),
\mcal{M}) . \]
Because each sheaf
$\mcal{H}om_{\mcal{O}_{X}} (\mcal{C}^{q}(X), \mcal{J}^{p})$
is flasque it follows that
\[ \mrm{H}^{q} \Gamma(X, \mcal{H}om_{\mcal{O}_{X}}
(\mcal{C}^{\bdot}(X), \mcal{J}^{\bdot})) =
\mrm{H}^{q} \mrm{R} \Gamma(X, \mrm{R} \mcal{H}om_{\mcal{O}_{X}}
(\mcal{C}^{\bdot}(X), \mcal{M})) . \]
The left hand side is by definition
$\mrm{H} \mrm{H}^{q}(X, \mcal{M})$. The right hand side is,
according to the theorem,
\[ \mrm{H}^{q} \mrm{R} \Gamma(X, \mrm{R}
\mcal{H}om_{\mcal{O}_{X^{2}}}(\mcal{O}_{X}, \mcal{M})) \cong
\opn{Ext}^{q}_{\mcal{O}_{X^{2}}}(\mcal{O}_{X}, \mcal{M})) . \]
\end{proof}

\section{Decomposition in Characteristic $0$}
\label{sec4}

In is section we prove that the Hochschild cochain complex
decomposes when $X$ is smooth and $\opn{char} \mbb{K} = 0$.
Throughout this section the base ring $\mbb{K}$ is
noetherian and $X$ is a separated finite type scheme over
$\mbb{K}$.

Let $A$ be a finitely generated $\mbb{K}$-algebra and
$\Omega^{q}_{A} = \Omega^{q}_{A / \mbb{K}}$
the module of relative K\"{a}hler differentials of degree $q$.
We declare $\boplus_{q} \Omega^{q}_{A}[q]$ to be a complex
with trivial coboundary.
For any $q \geq 0$ there is an $A$-linear homomorphism
\[ \pi: \mcal{C}_{q}(A) =
\mcal{B}_{q}(A) \otimes_{A^{\mrm{e}}} A \ar
\Omega^{q}_{A} , \]
\[ \pi((1 \otimes a_{1} \otimes \cdots \otimes a_{q} \otimes 1)
\otimes 1) = \mrm{d} a_{1} \wedge \cdots \wedge \mrm{d} a_{q} . \]
Since $\pi \partial = 0$ we obtain a homomorphism of complexes
$\pi: \mcal{C}^{\bdot}(A) \ar \boplus_{q} \Omega^{q}_{A}[q]$.

Recall that
$I_{q} = \opn{Ker}(\mcal{B}_{q}(A) \to A)$.

\begin{lem} \label{lem3.1}
Let $m > q$. Then
$\pi(I_{q}^{m} \cdot \mcal{C}_{q}(A)) = 0$.
\end{lem}

\begin{proof}
Let us consider $\Omega^{q}_{A}$ as a $\mcal{B}_{q}(A)$-module.
Then $\pi$ is a differential operator of order $\leq q$. Now use
\cite{Ye1} Proposition 1.4.6.
\end{proof}

The lemma shows that $\pi$ is continuous, so it extends to a
homomorphism of complexes
\[ \pi: \what{\mcal{C}}^{\bdot}(A) \ar
\boplus_{q} \Omega^{q}_{A}[q] . \]

If we take
$A = \mbb{K}[\bsym{t}] := \mbb{K}[t_{1}, \ldots, t_{n}]$
the polynomial algebra in $n$ variables, then
$\mcal{B}_{q}(\mbb{K}[\bsym{t}])$ is a polynomial algebra over
$\mbb{K}[\bsym{t}]$ in the $n (q + 1)$ elements
$\tilde{\mrm{d}}_{j} t_{i}$, cf.\
Lemma \ref{lem2.5}. Put a $\mbb{Z}$-grading on
$\mcal{B}_{q}(\mbb{K}[\bsym{t}])$ by declaring
$\opn{deg} \tilde{\mrm{d}}_{j} t_{i} := 1$, and
$\opn{deg} a := 0$ for $0 \neq a \in \mbb{K}[\bsym{t}]$.
This induces a grading on
$\mcal{C}_{q}(\mbb{K}[\bsym{t}]) =
\mcal{B}_{q}(\mbb{K}[\bsym{t}])
\otimes_{\mcal{B}_{0}(\mbb{K}[\bsym{t}])} \mbb{K}[\bsym{t}]$.
Also consider $\Omega^{q}_{\mbb{K}[\bsym{t}]}$ to be homogeneous
of degree $q$.

\begin{lem} \label{lem3.2}
The homomorphism
$\pi: \mcal{C}_{q}(\mbb{K}[\bsym{t}]) \to
\Omega^{q}_{\mbb{K}[\bsym{t}]}$
has degree $0$.
\end{lem}

\begin{proof}
Since $\mcal{C}_{q}(\mbb{K}[\bsym{t}])$ is a free
$\mbb{K}[\bsym{t}]$-module with basis the monomials
$\beta =$ \linebreak
$\tilde{\mrm{d}}_{j_{1}} t_{i_{1}} \cdots
\tilde{\mrm{d}}_{j_{m}} t_{i_{m}}$
with $1 \leq j_{1}, \ldots, j_{m} \leq q$
and $1 \leq i_{1}, \ldots, i_{m} \leq n$
it suffices to look at $\pi(\beta)$. We note that
$\opn{deg} \beta = m$. Now
\[ \pi((1 \otimes a_{1} \otimes \cdots \otimes a_{q} \otimes 1)
\otimes 1) =  0 \]
if any $a_{p} \in \mbb{K}$, $1 \leq p \leq q$. Therefore
$\pi(\beta) = 0$ unless
$\{ j_{1}, \ldots, j_{m} \} = \{ 1, \ldots, q \}$.
We conclude that $\pi(\beta) = 0$ if $m < q$. On the other hand,
since each $\tilde{\mrm{d}}_{j} t_{i} \in I_{q}$,
Lemma \ref{lem3.1} tells us that $\pi(\beta) = 0$ if $m > q$.
\end{proof}

The lemma says that $\pi$ is a morphism in the category
$\cat{GrMod} \mbb{K}[\bsym{t}]$ of $\mbb{Z}$-graded
$\mbb{K}[\bsym{t}]$-modules and degree $0$ homomorphisms.

\begin{lem} \label{lem3.3}
Assume $n!$ is invertible in $\mbb{K}$. Then
$\pi: \mcal{C}^{\bdot}(\mbb{K}[\bsym{t}]) \ar
\boplus_{q} \Omega^{q}_{\mbb{K}[\bsym{t}]}[q]$
is a homotopy equivalence of complexes over
$\cat{GrMod} \mbb{K}[\bsym{t}]$.
\end{lem}

\begin{proof}
Write $A := \mbb{K}[\bsym{t}]$. For $q > n$ we have
$\Omega^{q}_{A} = 0$, and $q!$ is invertible for all $q \leq n$.
So by \cite{Lo} Proposition 1.3.16 the
homomorphism of complexes
$\pi: \mcal{C}^{\bdot}(A) \ar \boplus_{q} \Omega^{q}_{A}[q]$
is a quasi-isomorphism. Now the complexes
$\mcal{C}^{\bdot}(A)$ and $\boplus_{q} \Omega^{q}_{A}[q]$
are both bounded above complexes of projective objects in
$\cat{GrMod} A$. So the quasi-isomorphism
$\pi: \mcal{C}^{\bdot}(A) \ar \boplus_{q} \Omega^{q}_{A}[q]$
has to be a homotopy equivalence. Namely there are homomorphisms
$\phi: \Omega^{q}_{A} \to \mcal{C}^{-q}(A)$
and
$h: \mcal{C}^{-q}(A) \to \mcal{C}^{-q - 1}(A)$
in $\cat{GrMod} A$ satisfying: $\partial \phi = 0$,
$1_{\mcal{C}^{-q}(A)} - \phi \pi = h \partial - \partial h$ and
$1_{\Omega^{q}_{A}} - \pi \phi = 0$.
\end{proof}

\begin{prop} \label{prop3.5}
Suppose $\mbb{K}[\bsym{t}] \ar A$ is \'etale and $n!$
is invertible in $\mbb{K}$. Then
$\pi: \what{\mcal{C}}^{\bdot}(A) \ar \boplus_{q}
\Omega^{q}_{A}[q]$
is a homotopy equivalence of topological $A$-modules. Namely there
are continuous $A$-linear homomorphisms
$\phi: \Omega^{q}_{A} \to \what{\mcal{C}}^{-q}(A)$
and
$h: \what{\mcal{C}}^{-q}(A) \to \what{\mcal{C}}^{-q - 1}(A)$
satisfying: $\partial \phi = 0$,
$1_{\what{\mcal{C}}^{-q}(A)} - \phi \pi = h \partial - \partial h$
and $1_{\Omega^{q}_{A}} - \pi \phi = 0$. Furthermore the
homomorphisms $\phi$ and $h$ are functorial in $A$.
\end{prop}

\begin{proof}
Declare $A$ to be homogeneous of degree $0$. From Lemma \ref{lem3.3}
we get homomorphisms
\[ \phi: A \otimes_{\mbb{K}[\bsym{t}]}
\Omega^{q}_{\mbb{K}[\bsym{t}]} \ar
A \otimes_{\mbb{K}[\bsym{t}]} \mcal{C}^{-q}(\mbb{K}[\bsym{t}]) \]
and
\[ h: A \otimes_{\mbb{K}[\bsym{t}]} \mcal{C}^{-q}(\mbb{K}[\bsym{t}])
\ar A \otimes_{\mbb{K}[\bsym{t}]}
\mcal{C}^{-q - 1}(\mbb{K}[\bsym{t}]) \]
in $\cat{GrMod} A$, satisfying the homotopy equations. Because
$\mbb{K}[\bsym{t}] \ar A$ is \'etale there is an isomorphism
$A \otimes_{\mbb{K}[\bsym{t}]} \Omega^{q}_{\mbb{K}[\bsym{t}]}
\cong \Omega^{q}_{A}$.
By Lemma \ref{lem2.5}, $\what{\mcal{C}}_{q}(A)$ is the
completion of
$A \otimes_{\mbb{K}[\bsym{t}]} \mcal{C}_{q}(\mbb{K}[\bsym{t}])$
with respect to the grading. Therefore $\phi$ and
$h$ extend uniquely to continuous homomorphisms as claimed. The
functoriality in $A$ follows from the uniqueness.
\end{proof}

\begin{thm} \label{thm3.1}
Let $\mbb{K}$ be a noetherian ring, let $X$ be a separated smooth
$\mbb{K}$-scheme of relative dimension $n$, and assume $n!$ is
invertible in $\mbb{K}$. Then for any complex
$\mcal{M}^{\bdot} \in \msf{D}(\cat{Mod} \mcal{O}_{X})$
the homomorphism of complexes
\[ \mcal{H}om_{\mcal{O}_{X}} \bigl(\boplus_{q} \Omega^{q}_{X}[q],
\mcal{M}^{\bdot} \bigr) \to
\mcal{H}om^{\mrm{cont}}_{\mcal{O}_{X}}
\bigl(\what{\mcal{C}}^{\bdot}(X), \mcal{M}^{\bdot}\bigr) \]
induced by $\pi$ is a quasi-isomorphism.
\end{thm}

\begin{proof}
The assertion may be checked locally on $X$, so let
$U = \opn{Spec} A \subset X$
be an affine open set admitting an \'{e}tale morphism
$U \to \mbf{A}^{n}_{\mbb{K}}$. If $U' = \opn{Spec} A' \subset U$
is any affine open subset then the ring homomorphisms
$\mbb{K}[\bsym{t}] \ar A \ar A'$ are \'etale. We deduce from
Proposition \ref{prop3.5} that
$\pi: \what{\mcal{C}}^{\bdot}(U) \ar
\boplus_{q} \Omega^{q}_{U}[q]$
is a homotopy equivalence of topological $\mcal{O}_{U}$-modules,
i.e.\ there are continuous $\mcal{O}_{U}$-linear homomorphisms
$\phi: \Omega^{q}_{U} \to \what{\mcal{C}}^{-q}(U)$
and
$h: \what{\mcal{C}}^{-q}(U) \to \what{\mcal{C}}^{-q - 1}(U)$
satisfying the homotopy equations.
\end{proof}

\begin{cor} \label{cor3.1}
Under the assumptions of the theorem, for any complex
$\mcal{M}^{\bdot} \in \msf{D}^{+}(\cat{Mod} \mcal{O}_{X})$
there is an isomorphism
\[ \boplus_{q} \bigl(\bwedge^{q} \mcal{T}_{X}\bigr)[-q]
\otimes_{\mcal{O}_{X}} \mcal{M}^{\bdot} \cong
\mrm{R} \mcal{H}om_{\mcal{O}_{X^{2}}}
(\mcal{O}_{X}, \mcal{M}^{\bdot}) \]
in $\msf{D}(\cat{Mod} \mcal{O}_{X^{2}})$. This isomorphism is
functorial in $\mcal{M}^{\bdot}$. In particular for
$\mcal{M}^{\bdot} = \mcal{O}_{X}$ we obtain
\[ \boplus_{q} \bigl(\bwedge^{q} \mcal{T}_{X}\bigr)[-q] \cong
\mrm{R} \mcal{H}om_{\mcal{O}_{X^{2}}}
(\mcal{O}_{X}, \mcal{O}_{X}) \]
in $\msf{D}(\cat{Mod} \mcal{O}_{X^{2}})$.
\end{cor}

\begin{proof}
Use Theorem \ref{thm2.1}.
\end{proof}

Observe that the isomorphism
$\bwedge^{q} \mcal{T}_{X} \cong
\mcal{E}xt_{\mcal{O}_{X^{2}}}^{q}(\mcal{O}_{X}, \mcal{O}_{X})$
deduced from Corollary \ref{cor3.1} differs by a factor of $q!$
from the Hochschild-Kostant-Rosenberg isomorphism (cf.\
\cite{HKR} Theorem 5.2 and \cite{Lo} Theorem 3.4.4).

Taking global cohomology in Corollary \ref{cor3.1} we deduce the next
corollary.

\begin{cor} \label{cor3.2}
Under the assumptions of the theorem, for any
$\mcal{O}_{X}$-module $\mcal{M}$ there is an isomorphism
\[ \opn{Ext}^{i}_{\mcal{O}_{X^{2}}}(\mcal{O}_{X}, \mcal{M})
\cong \boplus_{q} \mrm{H}^{i - q} \bigl( X,
\bigl(\bwedge^{q} \mcal{T}_{X}\bigr) \otimes_{\mcal{O}_{X}} \mcal{M}
\bigr) . \]
\end{cor}

Corollary \ref{cor3.2} was proved by Swan (\cite{Sw} Corollary 2.6)
in the case $X$ is smooth quasi-projective over the field
$\mbb{K} = \mbb{C}$.

Let us concentrate on the Hochschild cochain complex with values
in $\mcal{O}_{X}$. Here we give notation to the homomorphism induced
by $\pi$; it is

\[ \pi_{\mrm{cd}}: \bwedge^{q} \mcal{T}_{X} \to
\mcal{C}^{q}_{\mrm{cd}}(X) . \]
The precise formula on an affine open set $U = \opn{Spec} A$ is
\[ \pi_{\mrm{cd}}(v_{1} \wedge \cdots \wedge v_{q})
(1 \otimes a_{1} \otimes \cdots \otimes a_{q} \otimes 1) =
\sum_{\sigma \in \Sigma_{q}} \opn{sgn}(\sigma)
v_{\sigma(1)}(a_{1}) \cdots v_{\sigma(q)}(a_{q}) \]
for $v_{i} \in \mcal{T}_{A} = \opn{Der}_{\mbb{K}}(A)$ and
$a_{i} \in A$, where $\opn{sgn}(\sigma)$ denotes the sign
of the permutation $\sigma$.

Theorem \ref{thm3.1} says that $\pi_{\mrm{cd}}$ is a
quasi-isomorphism if $X$ is smooth of relative dimension $n$ and
$n!$ is invertible in $\mbb{K}$. The next result is a converse.

\begin{thm} \label{thm3.2}
Let $\mbb{K}$ be a Gorenstein noetherian ring of finite Krull
dimension and let $X$ be a smooth separated $\mbb{K}$-scheme of
relative dimension $n$. Then the following three conditions are
equivalent.
\begin{enumerate}
\rmitem{i} $\pi: \what{\mcal{C}}^{\bdot}(X) \ar
\boplus_{q} \Omega^{q}_{X}[q]$
is a quasi-isomorphism.
\rmitem{ii}
$\pi_{\mrm{cd}}: \boplus_{q} (\bwedge^{q} \mcal{T}_{X})[-q] \to
\mcal{C}^{\bdot}_{\mrm{cd}}(X)$
is a quasi-isomorphism.
\rmitem{iii} $n!$ is invertible in $\mcal{O}_{X}$.
\end{enumerate}
\end{thm}

\begin{proof}
All three conditions can be checked locally. So take a sufficiently
small affine open set $U = \opn{Spec} A \subset X$
such that there is an \'{e}tale homomorphism
$\mbb{K}[t_{1}, \ldots, t_{n}] \ar A$.
We will prove that the three conditions are equivalent on $U$
(cf.\ Propositions \ref{prop2.1} and \ref{prop2.2}).

\medskip \noindent (i) $\Leftrightarrow$ (ii):
Denote by $D$ the functor $\opn{Hom}_{A}(-, A)$
and by
$\mrm{R} D : \msf{D}(\cat{Mod} A) \ar \msf{D}(\cat{Mod} A)$
its derived functor. Consider the homomorphism of complexes
\[ \pi_{\mrm{cd}}: \boplus_{q}
\bigl(\bwedge^{q} \mcal{T}_{A}\bigr)[-q]
\ar \mcal{C}^{\bdot}_{\mrm{cd}}(A) . \]
By Lemma \ref{lem2.5}, $\what{\mcal{C}}_{q}(A)$ is a power series
algebra over $A$ in $n q$ elements. Hence the adjunction map
\[ \what{\mcal{C}}_{q}(A) \ar
D \mcal{C}^{q}_{\mrm{cd}}(A) =
\opn{Hom}_{A}(\opn{Hom}_{A}^{\mrm{cont}}(\what{\mcal{C}}_{q}(A),
A), A) \]
is bijective, and we get
\[ \pi = D(\pi_{\mrm{cd}}) : \what{\mcal{C}}_{q}(A) \ar
\Omega^{q}_{A} . \]
We claim that moreover
$\what{\mcal{C}}^{\bdot}(A) =
\mrm{R} D \mcal{C}^{\bdot}_{\mrm{cd}}(A)$
and
\begin{equation}
\pi = \mrm{R} D(\pi_{\mrm{cd}}) : \what{\mcal{C}}^{\bdot}(A) \ar
\boplus_{q} \Omega^{q}_{A}[q] .
\end{equation}
To verify this let us choose a bounded injective resolution
$A \ar J^{\bdot}$ in $\cat{Mod} A$, which is possible since $A$
is Gorenstein of finite Krull dimension. Each $A$-module
$\mcal{C}^{q}_{\mrm{cd}}(A)$ is free. Then, even though the
complex $\mcal{C}^{\bdot}_{\mrm{cd}}(A)$ is unbounded,
\[ \opn{Hom}_{A}(\mcal{C}^{\bdot}_{\mrm{cd}}(A), A) \ar
\opn{Hom}_{A}(\mcal{C}^{\bdot}_{\mrm{cd}}(A), J^{\bdot}) \]
is a quasi-isomorphism. Thus the claim is proved.

The functor $\mrm{R} D$ is a duality of the subcategory
$\msf{D}_{\mrm{c}}(\cat{Mod} A)$ of complexes with finitely
generated cohomologies. By Corollary \ref{cor2.1} we know that
$\mcal{C}^{\bdot}_{\mrm{cd}}(A) \in$ \linebreak
$\msf{D}_{\mrm{c}}(\cat{Mod} A)$, and clearly
$\boplus_{q} (\bwedge^{q} \mcal{T}_{A})[-q] \in
\msf{D}_{\mrm{c}}(\cat{Mod} A)$.
We conclude that $\pi_{\mrm{cd}}$ is an isomorphism in
$\msf{D}_{\mrm{c}}(\cat{Mod} A)$ iff
$\pi = \mrm{R} D(\pi_{\mrm{cd}})$ is an isomorphism.

\medskip \noindent (i) $\Leftrightarrow$ (iii):
We know that
$\mcal{C}^{\bdot}(A) \ar \what{\mcal{C}}^{\bdot}(A)$
is a quasi-isomorphism. Let
$\epsilon: \Omega^{q}_{A} \ar \mrm{H}^{-q} \mcal{C}^{\bdot}(A)$
be the isomorphism of the Hochschild-Kostant-Rosenberg Theorem
\cite{Lo} Theorem 3.4.4. Then by \cite{Lo} Proposition 1.3.16,
$\pi \epsilon(\alpha) = q! \alpha$ for all
$\alpha \in \Omega^{q}_{A}$.

If $n!$ is invertible in $A$ then so is $q!$ for all
$q \leq n$. For $q > n$ we have $\Omega^{q}_{A} = 0$. So $\pi$ is
a quasi-isomorphism.

Conversely, suppose $\pi$ is a quasi-isomorphism. Let $\alpha$
be a basis of the free $A$-module $\Omega^{n}_{A}$. Then
$n! \alpha = \pi \epsilon(\alpha)$ is also a basis,
so $n!$ must be invertible in $A$.
\end{proof}

Oddly, if $X$ is affine there is always a decomposition,
regardless of characteristic.

\begin{prop} \label{prop3.3}
If $\mbb{K}$ is noetherian and $X$ is affine and smooth over
$\mbb{K}$ then there is a canonical isomorphism
\[ \mrm{R} \mcal{H}om_{\mcal{O}_{X^{2}}}
(\mcal{O}_{X}, \mcal{O}_{X}) \cong
\boplus_{q} \bigl(\bwedge^{q} \mcal{T}_{X}\bigr)[-q] \]
in $\msf{D}(\cat{Mod} \mcal{O}_{X^{2}})$.
\end{prop}

\begin{proof}
Say $X = \opn{Spec} A$. Let $A \ar J^{\bdot}$ be an injective resolution
in $\cat{Mod} A^{\mrm{e}}$, and set
$N^{\bdot} := \opn{Hom}_{A^{\mrm{e}}}(A, J^{\bdot})$,
which is a complex of $A$-modules. Denote by
$F: \cat{Mod} A \to \cat{Mod} A^{\mrm{e}}$
the restriction of scalars functor for the homomorphism
$A^{\mrm{e}} \to A$ (this is the ring version of $\Delta_{*}$).
Then
$F N^{\bdot} = \mrm{R} \opn{Hom}_{A^{\mrm{e}}}(A, A)$
in $\msf{D}(\cat{Mod} A^{\mrm{e}})$.
Let
$G: \cat{Mod} A^{\mrm{e}} \to \cat{Mod} \mcal{O}_{X^{2}}$
be the sheafication functor. Since $G J^{\bdot}$ is an injective
resolution of $\mcal{O}_{X}$ we see that
\[ G F N^{\bdot} \cong
\mcal{H}om_{\mcal{O}_{X^{2}}} (\mcal{O}_{X}, G J^{\bdot}) \cong
\mrm{R} \mcal{H}om_{\mcal{O}_{X^{2}}} (\mcal{O}_{X}, \mcal{O}_{X})
\]
in $\msf{D}(\cat{Mod} \mcal{O}_{X^{2}})$.

Now according to the Hochschild-Kostant-Rosenberg Theorem (see
\cite{HKR} Theorem 5.2 and \cite{Lo} Theorem 3.4.4) the cohomology
$\mrm{H}^{q} N^{\bdot} = \opn{Ext}^{q}_{A^{\mrm{e}}}(A, A) \cong
\bwedge^{q} \mcal{T}_{A}$.
Since the $A$-modules
$\bwedge^{q} \mcal{T}_{A}$ are projective and almost all of them
are zero, it is easy to see, by truncation and splitting, that
$N^{\bdot} \cong \boplus_{q} (\bwedge^{q} \mcal{T}_{A})[-q]$
in $\msf{D}(\cat{Mod} A)$. Therefore
$G F N^{\bdot} \cong
\boplus_{q} \bigl(\bwedge^{q} \mcal{T}_{X}\bigr)[-q]$
in $\msf{D}(\cat{Mod} \mcal{O}_{X^{2}})$.
\end{proof}

\begin{que}
We have seen that if $X$ is affine or if $\mbb{K}$ contains enough
denominators then the Hochschild cochain complex
$\mcal{C}^{\bdot}_{\mrm{cd}}(X)$ decomposes in the derived
category. Is there decomposition in other circumstances?
\end{que}

\begin{que}
How is the decomposition of Theorem \ref{thm3.1}
related to the Hodge decomposition of Gerstenhaber-Schack
\cite{GS}?
Perhaps the comparison to Swan's definition of
Hochschild cochains (Section \ref{sec3}) can help.
\end{que}

\begin{rem}
In \cite{Ko},
$\mcal{D}^{\bdot}_{\mrm{poly}}(X) :=
\mcal{C}^{\bdot}_{\mrm{cd}}(X)[1]$
is called the complex of poly-differential operators. The complex
$\mcal{T}^{\bdot}_{\mrm{poly}}(X) :=
\boplus_{q} (\bwedge^{q} \mcal{T}_{X})[1 - q]$
is called the complex of poly-vector fields. Kontsevich's Formality
Theorem \cite{Ko} says that
$(\frac{1}{q!} \pi^{q}_{\mrm{cd}})_{q \geq 0}$
is the degree $1$ component of an
$\mrm{L}_{\infty}$-quasi-isomorphism of
the DG Lie algebra structures of $\mcal{D}^{\bdot}_{\mrm{poly}}(X)$
and $\mcal{T}^{\bdot}_{\mrm{poly}}(X)$
when $\mbb{K}$ is a field of characteristic $0$.
\end{rem}

\end{document}